\documentclass[12pt]{article}

\usepackage{amssymb,amsmath,amsfonts,amssymb}
\usepackage{graphics,graphicx,color}

\def\@abssec#1{\vspace{.05in}\footnotesize \parindent .2in
{\bf #1. }\ignorespaces}

\graphicspath{{/EPSF/}{Figures/}}
\DeclareGraphicsExtensions{.eps}

\def \Rm {\mathbb R}

\def \Zm {\mathbb Z}
\def \Sm {\mathbb S}

\newcommand{\pdr}[2]{\dfrac{\partial{#1}}{\partial{#2}}}

\newcommand{\bx}{\mathbf x}

\newcommand{\mM}{\mathcal M}

\newcommand{\cout}[1]{}

 \renewcommand{\arraystretch}{1.5}

\title{Reconstruction of Cloud Geometry from High-Resolution Multi-Angle 
Images}

\author{Guillaume Bal \thanks{Department of Applied Physics and
        Applied Mathematics, Columbia University,
        New York NY 10027, USA; gb2030@columbia.edu} \and JiaMing Chen \thanks{Department of Mathematical Sciences, Rensselear Polytechnic Institute, Troy, NY 12180, USA} \and Anthony B. Davis \thanks{Jet Propulsion Laboratory, California Institute of Technology, Pasadena, CA
91109, USA}}

\begin{document}

\maketitle


\begin{abstract}
  We consider the reconstruction of the interface of compact, connected ``clouds" from satellite or airborne light intensity measurements. In a two dimensional setting, the cloud is modeled by an interface, locally represented as a graph, and an outgoing radiation intensity that is consistent with a diffusion model for light propagation in the cloud. Light scattering inside the cloud and the optical internal parameters of the cloud are not modeled. The main objective is to understand what can or cannot be reconstructed in such a setting from intensity measurements in a finite (on the order of 10) number of directions along the path of a satellite. Numerical simulations illustrate the theoretical predictions.
\end{abstract}


\renewcommand{\thefootnote}{\fnsymbol{footnote}}
\renewcommand{\thefootnote}{\arabic{footnote}}

\renewcommand{\arraystretch}{1.1}





\section{Introduction}
\label{sec:intro}

This paper concerns the reconstruction of a cloud surface from radiance satellite measurements.
The state-of-the-art for extraction of physical properties of clouds from passive remote sensing is based on a pixel-by-pixel interpretation of the measured radiance across the solar spectrum, possibly at different angles and states of polarization.  Invariably, the cloud is modeled as a plane-parallel slab \cite{MK-JAS-90,Petal-GRS-03}.  This radical assumption about cloud geometry is more-or-less justified for stratiform clouds that have considerably more horizontal extension than thickness in the vertical.  As important as these clouds are for the balance of the Earth's climate, that leaves out many important types of clouds resulting typically from more vigorous convection (updrafts)---the cumulus class of cloud types.  For such clouds, the first order of business in remote sensing should be to reconstruct their non-trivial outer shape, hence the goal of the present demonstration.  The main objective of the paper is to reconstruct such an interface as well as the light intensity emitted from it without modeling the internal optical properties of the cloud.

How light propagates inside the cloud is accurately modeled by a radiative transfer model \cite{chandra,MD-SP-05,DM-RPP-10}.  The plane-parallel slab geometry is attractive largely because it leads to one-dimensional radiative transfer that is considered to be a tractable problem in computational physics \cite{chandra}.  As a demonstration that one can depart radically from the plan-parallel cloud geometry and still have an analytically tractable radiative transfer problem to solve, at least in the relevant diffusion approximation, Davis \cite{Davis2002SPIE} assumed perfectly spherical clouds; he then used his solution to derive effective optical thicknesses---or rather, diameters---for cumulus in a sparse field of broken clouds.
Reconstruction of the constitutive (optical) parameters in a full three-dimensional radiative transfer equation is a notoriously difficult and ill-posed problem \cite{B-IP-09} that has rarely been attempted for real clouds. For computationally intense efforts in that direction during the past decade or so, see \cite{MarchandEtAl2004JGR,CornetEtAl2008JGR}. Far more efficient methods that may lead to practical implementations are currently being investigated \cite{LSAD-15,MCB-JQSRT-14}.

In this paper, we do not aim to reconstruct such parameters and rather identify geometric properties of the cloud that are directly observable from satellite measurements and are independent from its optical properties.  We restrict ourselves to a two-dimensional setting to simplify the modeling and the numerical simulations. We expect most results presented here to hold without major modification in a three dimensional environment. The cloud is modeled as a compact, connected, domain with a sufficiently smooth interface. The interface of the cloud is  given by a curve parameterized by $t\mapsto\gamma(t)\in\Rm^2$. Since the reconstructions are local, the curve is also represented locally by a graph $(x,h(x))$. The main objective is then to estimate $h(x)$ from available data. To first order, the radiative transfer of light inside a cloud goes as follows. Light is emitted from the sun, enters the cloud proper, scatters (multiple times), and leaves the cloud at the points $(x,h(x))$ in all possible outgoing directions $\theta$. Importantly for the cloud shape reconstruction, cloud-escaping radiation is not isotropic (a.k.a. Lambertian).

The measurements obtained by sensors mounted on satellites are modeled as the light intensity at position $(x,Z)$ for all $x$ and a fixed $Z$ (the satellite's orbital altitude above the Earth's surface). We assume that the sensors are equipped with several directional filters so that light intensity $u(x,Z,\theta)$ is measured for a discrete number of directions $\theta\in\{\theta_j\}_{1\leq j\leq J}$. Typically, $J=9$, as has been implemented for the Multi-Angle Imaging Spectro-Radiometer (MISR) \cite{Detal-GRS-98}.

MISR is part of the payload on NASA's Terra satellite, which has a morning-crossing sun-synchronous orbit at $Z$ = 705~km.  MISR's nine imaging sensors are ``push-broom'' cameras that use the orbital motion pitched at fixed angles ranging from about 70$^\circ$ in the forward direction to the same in the backward direction.  Their common spatial sampling is 0.275~km.  Thus a typical cumulus cloud, which has commensurate horizontal and vertical dimensions of of a few km, will cover a few tens of pixels in the along-track direction.  In the present study, we will assume many more spatial samples.  This is more like what can be achieved using a digital camera located on the ground while the cloud is advected past it by a steady wind ($Z$ is then related to cloud height), or else by using an airborne platform ($Z$ is then related to the aircraft's altitude).  There are space-based imagers that can achieve very high ($\sim$10~m) resolution, but they only have one view angle.  Therefore, only quite recently developed airborne imaging sensors have multi-angle capability \cite{Detal-AMT-13}.  We will address the issues raised by the limited spatial sampling in current space-based imaging sensors in future work and, for the moment, we will continue to call the platform a ``satellite'' since satellites will eventually be the source of abundant free data for cloud shape reconstruction.

Let $\nu(x)$ be the outward unit normal vector at $(x,h(x))$, which is thus given by the proper normalization (by $(1+(h'(x))^2)^{-\frac12}$) of the vector $(-h'(x),1)$.  Here $h'$ is the derivative of $h$. The outward radiation at a point $(x,h(x))$ is therefore characterized by the intensity $u(x,h(x),\theta)$ for all vectors $\theta$ such that $\theta\cdot\nu(x)>0$. This is a two-dimensional function of $(x,\theta)$ whereas the available measurements correspond to $J$ one-dimensional functions. This results in a severely under-determined problem and assumptions on the outgoing radiation intensity are necessary. In this paper, we assume that the outgoing light intensity is of the form $u(x,h(x),\theta)=\alpha(x)H(\theta\cdot\nu(x))$ where $\alpha$ and $H$ are therefore functions of one variable. This physics-based model for light intensity is the correct one when light propagation in the cloud is modeled as a diffusive process \cite{chandra,B-IP-09}. This situation holds for sufficiently opaque clouds, which are the ones for which a sharp separation between the inside and outside of the cloud is the most realistic.

We are fully aware, as hypothesized in \cite{Mandelbrot1977} and proven observationally in \cite{Lovejoy1982} that, due to the turbulent nature of their dynamics, cloud shapes are best represented as fractals, that is, as convoluted surfaces that do not have well-defined tangent planes nor normal vectors at any scale.  Consequently, our methodology proposes to deliver a smoothed or better yet, properly averaged, approximation of the cloud's actual outer three-dimensional shape.  Such smoothing is necessary in the context of limited available data and remains a vast improvement over the aforementioned operational assumption that clouds are plane-parallel optical media, regardless of image context.

The main result of this paper is an iterative reconstruction procedure for $(h,\alpha,H)$ when $J$ is sufficiently large based on a linearization of the functional mapping $(h,\alpha,H)$ to the available satellite data. In ``favorable'' situations, which depend on the state about which the linearization is performed, this linear map is invertible and provides stable reconstructions for $(h,\alpha,H)$. Explicit calculations allow us to display several such ``favorable'' as well as less-favorable situations. Numerical simulations confirm the theoretical predictions.

The outline for the rest of the paper is as follows.
The graph model for the cloud interface and the associated inverse problems are presented in section \ref{sec:graph}. The numerical algorithm used to solve the inverse problem is given in section \ref{sec:num}.  In section \ref{sec:speed}, we show that the available measurements are not capable of uniquely reconstructing the parameters $(h,\alpha,H)$ when the (scalar-valued) speed of the cloud is also unknown. Finally, a more general geometric setting based on a polar representation of the cloud geometry is given in section \ref{sec:geom} along with numerical simulations.  A model of light scattering in the cloud, or different hypotheses on the structure of the outgoing radiation, are then necessary to uniquely reconstruct $(h,\alpha,H)$ as well as the cloud speed.

\section{Graph model for cloud reconstruction}
\label{sec:graph}

We present the geometric assumptions on the cloud in section \ref{sec:gs}. The inverse problem is then linearized in section \ref{sec:lin} and an iterative algorithm is described in section \ref{sec:nonlin}. The calculations of section \ref{sec:lin} are the main result of this paper and show under which conditions all or some of the parameters of interest $(h,\alpha,H)$ can be reconstructed.

\subsection{Geometric setting}
\label{sec:gs}

\begin{figure}\begin{center}
\includegraphics[height=8cm]{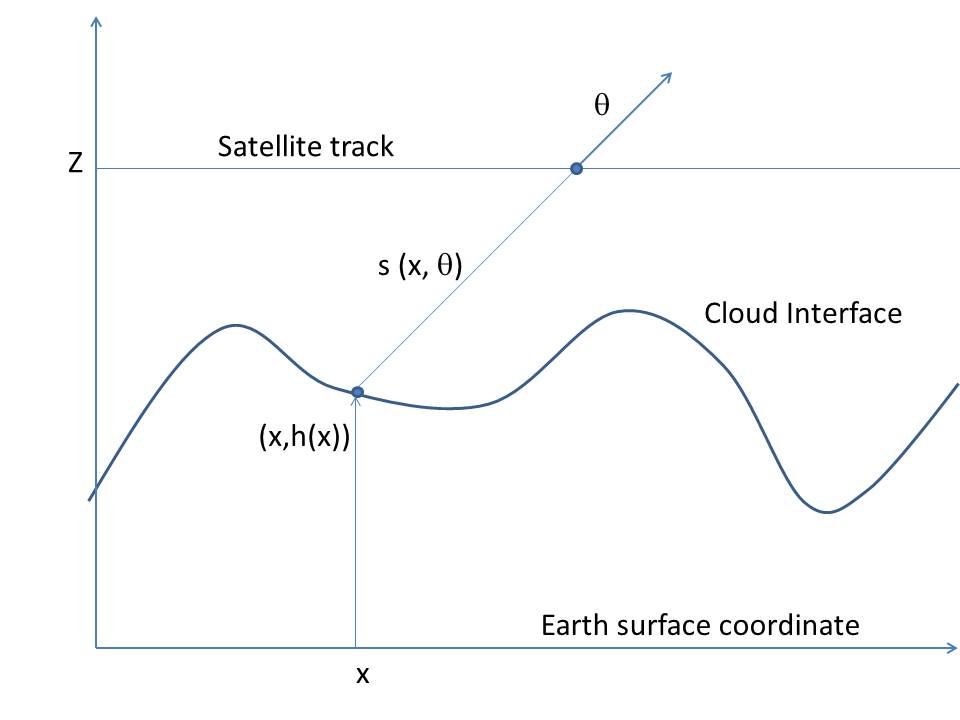}
\end{center}
\caption{Geometry of cloud interface}
\label{fig:1}
\end{figure}

We consider the setting of a two-dimensional bounded domain (the cloud) with boundary given by a curve parameterized by $t\mapsto\gamma(t)\in\Rm^2$. To simplify the presentation, we focus in this section on the reconstruction of the ``upper" part of the cloud, which we assume is given by the graph of a function $h(x)$; in other words,  $\gamma(x)=(x,h(x))\in\Rm^2$, with $h(x)$ the unknown function we wish to reconstruct. Since reconstructions will be shown to be local, an arbitrary boundary may be reconstructed by using the appropriate number of local graphs representing the interface; see also Fig.~\ref{fig:1}.

The available measurements are assumed to be obtained by an imaging sensor at a fixed elevation $Z$ (we neglect surface curvature, which can be accounted for in a straightforward manner). Measurements are performed along the line $X\mapsto (X,Z)$ and for a given number of (upward-pointing) directions $\theta_j=(\cos\phi_j,\sin\phi_j)$ for $\phi_j\in[0,\pi]$ and $1\leq j\leq J$. A typical practical value is $J=9$ \cite{Detal-GRS-98,Detal-AMT-13}.

Most techniques for the reconstruction of strongly scattering optical media 
are based on the assumption that the equations of radiative transfer are valid to model the propagation of photons inside the medium. These models require to first reconstruct spatially varying functions such as the diffusion coefficient (scattering mean free path) and the absorption coefficient. Such parameter estimations are often ill-posed. When satellite measurements are available for a large number of angles, we argue that the reconstruction of such parameters can be bypassed in some configurations, at least in a first stage. We propose an inversion procedure whose main objective is a direct reconstruction of $h(x)$. This requires some assumptions on the radiation emitted from the cloud that we now present in detail.

Let $u(x,z,\theta)$ be the (phase-space) density of photons as a function of position $(x,z)$ and angle $\theta$. We assume that the medium between the cloud and the satellite is clear, so that photons advect freely:
\begin{align}
 \theta\cdot\nabla u(x,z,\theta) &=0 &&& \nabla &= (\pdr{}x, \pdr{}z )^t.
\end{align}
This implies that for all points $(X,Z)-s\theta$ outside of the cloud, we have
\begin{displaymath}
   u(X,Z,\theta) = u(X-s\cos\phi,Z-s\sin\phi,\theta).
\end{displaymath}
For a point $(X,Z)$ on the satellite line of measurements and for $\theta$ pointing ``upwards", we denote by $s(X,Z,\theta)$ the distance from the satellite to the cloud in the direction $\theta$ and by $(x(X,Z,\theta),h(x(X,Z,\theta)))$ the point of intersection on the surface of the cloud. These functions are defined via the constraints
\begin{align}\label{eq:sx}
X-s(X,Z,\theta)\cos\phi &= x(X,Z,\theta), & Z-s(X,Z,\theta)\sin\phi &= h(x(X,Z,\theta)).
\end{align}
In our theoretical analyses, we assume that the functions $s$ and $x$ are uniquely defined and smooth and we set $s=+\infty$ if the half-line from $(X,Z)$ in direction $\theta$ does not intersect the cloud. For a fixed $\theta$, we thus assume that the mapping $X\to x(X,Z,\theta)$ is a local diffeomorphism.

Information from the satellite measurements thus provides knowledge of
\begin{align}\label{eq:sat1}
   u(X,Z,\theta) = u(x(X,Z,\theta),h(x(X,Z,\theta),\theta).
\end{align}
The outgoing radiation at the cloud surface is given by $u(x,h(x),\theta)$. It is not possible to reconstruct an arbitrary radiation profile $u(x,h(x),\theta)$ from knowledge of $(X,Z,\theta)\mapsto u(X,Z,\theta)$ since in fact we find a solution $u(x,h(x),\theta)$ for each choice of $h(x)$. The outgoing radiation at the cloud surface thus needs to be constrained.

Assuming that scattering (disorder) is large in the cloud (i.e., the scattering mean free path is small compared to the overall extension of the cloud, as is often the case), then photon propagation inside the cloud is well-approximated by a diffusion equation.  In that case, $u(x,z,\theta) \sim U(x,z) + 2 \theta \cdot F(x,z)$ with $\{U,F\}$ being the solution of equations of the form $\nabla \cdot F + \sigma_a U = 0$ and $F = -D \nabla U$.  This approximation breaks down at the surface of the cloud, where it needs to be replaced by the solution of a Milne problem \cite{chandra,B-IP-09}. At a point $(x,h(x))$ of the surface, let
\begin{displaymath}
   \nu(x) = \dfrac1{\sqrt{1+(h'(x))^2}} \,\big(-h'(x),1\big)
\end{displaymath}
be the outward unit normal vector to the surface. Then we find the model
\begin{align}
  u(x,h(x),\theta) = U(x,h(x)) H(\theta\cdot\nu(x)),
\end{align}
with $U$ the solution of the diffusion approximation and $H$ an appropriate Chandrasekhar function that depends on the scattering phase function inside the cloud \cite{chandra,B-IP-09}. We note that $\theta\cdot\nu=\sin(\phi-\mu)$ with $\tan\mu=h'$. Therefore, $H(\theta\cdot\nu(x))= H\circ \sin (\phi-\arctan h'(x))$.


Our constraint on the outgoing radiation of the cloud is to be consistent with the above diffusion approximation. In short, we assume that
\begin{align}\label{eq:chbeta}
 u(x,h(x),\theta) = \alpha(x) \beta(\phi-\arctan h'(x)),\qquad
 \beta := H\circ \sin,
\end{align}
for two unknown functions: $\alpha(x)$, a measure of the angularly integrated radiance escaping the cloud at point $x$; and $\beta(\theta):=H\circ \sin (\theta)$, an angular distribution model that is assumed independent of position in the sense that it is function only of the angle of the observed radiance with the local normal. From \eqref{eq:sat1} we thus observe that
\begin{align}
\label{eq:sat2}
 u(X,Z,\theta) =  \alpha \big(x(X,Z,\theta)\big) \beta \big(\phi-\arctan h'(x(X,Z,\theta))  \big),
\end{align}
with $x(X,Z,\theta)$ and $\nu (x(X,Z,\theta))$ functionals of the unknown height function $h(x)$.

Note that multiplying $\alpha$ by a constant and dividing $\beta$ by the same constant does not change $u(X,Z,\theta)$. We thus assume that $\beta(\phi)$ is normalized, for instance by assuming that it integrates to $1$ or that the value for a fixed $\phi$ is known.

In our setting, we thus have a measurement operator
\begin{align}
\label{eq:measop}
{\mathcal M}: \, (\alpha,\beta,h) \mapsto {\mathcal M}(\alpha,\beta,h)  = u(X,Z,\theta).
\end{align}
The inverse problem consists of reconstructing $(\alpha,\beta,h)$ from knowledge of ${\mathcal M}(\alpha,\beta,h)  = u(X,Z,\theta)$ for a fixed value of $Z$, for $X$ on a line segment, and for $\theta\in \{\theta_j\}_{1\leq j\leq J}$.


\subsection{Linearization of the inverse problem}
\label{sec:lin}

Let us define $v:=(\alpha,\beta,h)$. The above measurement operator ${\mathcal M}$, a nonlinear functional of $v$, does not seem to have an explicit inversion formula. Following a standard methodology, we linearize the problem about a reference configuration $v_0=(\alpha_0,\beta_0,h_0)$ and obtain an equation for the linear update $\delta v=(\delta \alpha,\delta \beta,\delta h)$ from knowledge of $\delta\mM = {\mathcal M}(\alpha,\beta,h) -{\mathcal M}(\alpha_0,\beta_0,h_0) =\delta u(X,Z,\theta)$.

\subsubsection{Linearized system of equations}

This section presents the main results of the paper, namely equation \eqref{eq:linear} providing the relationship between the unknown $\delta v$ and the known measurement $\delta u$.

Recall \eqref{eq:sat2}:
  $u(X,Z,\theta) =  \alpha \big(x(X,Z,\theta)\big) \beta \big(\phi-\arctan h'(x(X,Z,\theta))  \big).$
We wish to differentiate this expression with respect to $(\alpha,\beta,h)$. While the derivatives with respect to $\alpha$ and $\beta$ are straightforward, the differentiation with respect to $h$ requires a few additional steps.

Since the graph $(x_0,h_0(x_0))$ is assumed to be known in the linearization procedure, we define $s_0(X,Z,\theta)$ and  $x_0(X,Z,\theta)$ as the solutions to \eqref{eq:sx} with $h$ replaced by $h_0$. Let $\delta s$ and $\delta x$ be the linearizations of $s(X,Z,\theta)-s_0(X,Z,\theta)$ and $x(X,Z,\theta)-x_0(X,Z,\theta)$, respectively.   We use the notation $f(x_0)$ for $f(x_0(X,Z,\theta))$ and identify $f(\theta)=f(\phi)$ recalling that $\theta=(\cos\phi,\sin\phi)$.

Let us introduce $\psi_0(x)=\arctan h'_0(x)$. Differentiating the above expression for $u$, we find to leading order that
\begin{displaymath}
\begin{split}
 \delta u &= \beta_0\big(\phi-\psi_0(x_0)\big) \,\delta \alpha(x_0) + \alpha_0(x_0) \,\delta\beta(\phi-\psi_0(x_0)) \\ &
  + \alpha'_0(x_0)\beta_0(\phi-\psi_0(x_0))  \delta x(x_0) \\
  &+  \alpha_0(x_0)    \beta'_0(\phi-\psi_0(x_0)) \Big(
       \dfrac{h_0^{''}(x_0)}{1+(h_0'(x_0))^2} \delta x(x_0)
      + \dfrac{-1}{1+(h_0'(x_0))^2}  (\delta h)'(x_0) \Big).
\end{split}
\end{displaymath}

It remains to express $\delta x$ in terms of $\delta h$. From \eqref{eq:sx}, we find up to second-order that
\begin{displaymath}
  -\delta s(x_0,\phi) \cos\phi = \delta x(x_0,\phi),\quad -\delta s(x_0,\phi) \sin\phi = h'_0(x_0)\delta x (x_0,\phi) + \delta h(x_0).
\end{displaymath}
We thus deduce that
\begin{align}
  \delta x(x_0,\phi) = \dfrac{\delta h(x_0)}{\tan \phi - h'_0(x_0)}.
\end{align}
Then, from the above expressions for $\delta u$ and $\delta x$, we obtain that
\begin{align}
\label{eq:linear}
  \delta u(X,Z,\phi) & = \beta_0\big(\phi-\psi_0(x_0)\big) \,\delta \alpha(x_0) + \alpha_0(x_0) \,\delta\beta(\phi-\psi_0(x_0)) \\ & + \dfrac{\alpha'_0(x_0)\beta_0(\phi-\psi_0(x_0)) + \alpha_0(x_0)  \psi'_0(x_0)\beta'_0(\phi-\psi_0(x_0))}{\tan \phi - h'_0(x_0)} \,\delta h(x_0) \notag \\
  & +    \beta'_0(\phi-\psi_0(x_0)) \dfrac{-\alpha_0(x_0) }{1+(h_0'(x_0))^2}  (\delta h)'(x_0),\notag
\end{align}
recalling that $\psi_0(x)=\arctan h'_0(x)$. 

The above expression is the main calculation of the paper.  Our objective is to reconstruct two spatially dependent functions $\delta\alpha(x)$ and $\delta h(x)$ and one angularly dependent function $\delta\beta(\phi)$ from the spatially and angularly dependent function $\delta u(X,Z,\phi)$. Counting dimensions shows that the problem is formally over-determined as soon as measurements are available for $\{\phi_j\}_{1\leq j\leq J}$ for $J\geq3$. However, rewriting the above equation as $A\delta v =\delta u$, and in normal form
\begin{align}
\label{eq:normal}
A^t A \delta v = A^t \delta u
\end{align}
with $A^t$ the adjoint operator to $A$, the linear operator $A^tA$ needs to be invertible. This imposes some conditions on $v_0=(\alpha_0,\beta_0,h_0)$, which are not always met in practice.
The subsequent sections present situations in which we can prove that $A^tA$ is indeed invertible provided that adapted boundary conditions (for instance Dirichlet conditions for $h(x)$) are imposed.

\subsubsection{Analysis of the linearized system}

As recalled earlier, we assume that the map $X\to x_0(X,Z,\theta)$ for each fixed $\theta$ is a smooth change of variables.  We easily verify that such is the case when $h'_0(x)$ is sufficiently small. Then the map $(X,\phi)\to(x_0(X,\phi),\phi+\psi_0(x_0))$ is also a smooth invertible change of variables. We then identify functions $u(X,\phi)$ and $u(x_0,\phi)$ with $\phi\leftarrow\phi+\psi_0(x_0)$ and use this convenient change of notation.

Let us consider that $\delta u(X,Z,\phi)$ is available for a continuum of values of $\phi$ to simplify the analysis.
Dividing \eqref{eq:linear} by $\alpha_0(x_0)\beta_0(\phi)$, we may recast it as
\begin{align}
\label{eq:linearbis}
  \dfrac{\delta u}{\alpha_0\beta_0}(x_0,\phi) & = (\delta \ln\alpha)(x_0) + (\delta \ln\beta)(\phi)
  + \psi_3(x_0,\phi) \delta h(x_0)   + \psi_4(x_0,\phi) (\delta h )' (x_0),
\end{align}
with $\delta \ln\alpha \sim \frac{\delta \alpha}{\alpha_0}$ up to lower-order terms, and with $\psi_3$ and $\psi_4$ coming explicitly from \eqref{eq:linear}.  The reconstruction of $(\delta \ln\alpha,\delta \ln\beta,\delta h)$ depends on the independence of the functions $(1,\psi_3,\psi_4)$
appearing in the system \eqref{eq:linearbis}.

 Let us consider first the simplest setting with $\beta_0(\phi)=\sin\phi$ (as a simple model for limb-darkening) and $h'_0(x)=0$ (flat horizontal cloud boundary). We verify that
\begin{displaymath}
 \psi_3(x_0,\phi) = \dfrac{(\ln \alpha_0(x_0))'}{\tan \phi},\quad \psi_4(x_0,\phi) = -(\ln \beta_0(\phi))' = -\dfrac{1}{\tan\phi}.
\end{displaymath}
In other words, we have
\begin{displaymath}
\dfrac{\delta u}{\alpha_0\beta_0}(x_0,\phi) = (\delta \ln\alpha)(x_0) + (\delta \ln\beta)(\phi)  + \cot \phi \big( (\ln \alpha_0)'(x_0) \delta h - (\delta h )' \big).
\end{displaymath}
If we normalize $\beta(\frac\pi2)=1$, then we observe that the above expression at $\phi=\frac\pi2$ (nadir view) provides $(\delta \ln\alpha)(x_0)$, or equivalently $\delta\alpha(x_0)$. It is then straightforward to observe that $(\delta \ln\beta)(\phi)$ is known up to the addition of $\lambda\cot \phi$ for $\lambda\in\Rm$ arbitrary, and that $(\ln \alpha_0(x_0))' \delta h - (\delta h )' $ is known up to the addition of the constant $\lambda$.

Consider the problem of the reconstruction of the parameters on a domain $x_{\rm min}<x_0<x_{\rm max}$. Then all parameters are uniquely determined provided that, for instance, $h(x_{\rm min})$ and $h(x_{\rm max})$ are known. Indeed, in such a setting with the additional constraint $-(\delta h)'' + \big((\ln \alpha_0(x_0))' \delta h\big)'=S$ with $S$ a known source, $\delta h$ is uniquely reconstructed, and then so is $\delta\beta$.

Consider a more general case with still $h'_0=0$ but $\beta_0$ arbitrary so that
\begin{displaymath}
\dfrac{\delta u}{\alpha_0\beta_0}(x_0,\phi) = (\delta \ln\alpha)(x_0) + (\delta \ln\beta)(\phi)  + \cot\phi (\ln \alpha_0)'(x_0) \delta h    - (\ln\beta_0)'(\phi)  (\delta h )' .
\end{displaymath}
Let us assume that $(1,\cot\phi,(\ln\beta_0)'(\phi))$ are linearly independent functions. Then as before, $(\delta \ln\alpha)(x_0)+\lambda_1$, $(\ln \alpha_0)'(x_0) \delta h+\lambda_2$, and $(\delta h )'+\lambda_3$ are reconstructed with the constants $\lambda_j$ still unknown, while $(\delta \ln\beta)(\phi)+\lambda_1+\lambda_2\cot\phi-\lambda_3(\ln\beta_0)'(\phi)$ is also reconstructed.

In favorable situations, we need to impose less conditions than when $\beta_0(\phi)=\sin(\phi)$. Indeed, assume that $\delta h(x_{\min})=0$ (because $h(x_{\rm min})$ is known). If $(\ln \alpha_0)'(x_{\min})\not=0$, then $\lambda_2$ is known and we have an expression for $(\delta h )'(x_{\rm min})$. This provides an expression for $\lambda_3$. The normalization of $\delta\beta(\frac\pi2)=0$ then allows us to reconstruct $\lambda_1$.

The general case with $h'_0\not=0$ proceeds as above: we need to impose one or two conditions on $\delta h$ depending on whether $(1,\psi_3,\psi_4)$ are linearly independent or not. The above analysis shows that $(\delta\alpha,\delta\beta,\delta h)$ can be uniquely reconstructed with minimal additional information required on $\delta h$.

Let us conclude this section by considering the case where $\beta(\phi)$ is known, for instance chosen as $H\circ\sin$ with $H$ the Chandrasekhar function corresponding to the appropriate scattering phase function inside the cloud.  Then in the setting with $\cot\phi$ and $(\ln\beta)'(\phi)$ linearly dependent, we obtain a unique reconstruction of $\delta\alpha$ and $(\ln \alpha_0(x_0))' \delta h - (\delta h )'$. This provides a reconstruction of $\delta h$ if for instance $\delta h(x_{\min})=0$.  If $\cot\phi$ and $(\ln\beta)'(\phi)$ linearly independent, then $(\ln\alpha_0)'(x_0) \delta h$ and $(\delta h )'$ can both be reconstructed and $\delta h$ is then uniquely determined provided that $(\ln\alpha_0)'$ does not vanish uniformly. Thus in practice, we expect to uniquely reconstruct $(\delta\alpha,\delta h)$ from the available satellite measurements.

\subsection{Nonlinear inverse problem}
\label{sec:nonlin}

Let us assume that the nonlinear operator $A^tA$ as defined above is invertible. Then it will be invertible for nearby values of $v_0$ by continuity. We may then set up the following Newton inversion for the nonlinear problem. Let us assume that $v^{(k)}=(\alpha^{(k)},\beta^{(k)},h^{(k)})$ has been reconstructed and define $A^{(k)}$ as the linearization of $\mM$ in the vicinity of $v^{(k)}$. Then up to second-order we have
\begin{displaymath}
  \mM(v) - \mM(v^{(k)}) \sim A^{(k)} (v-v^{(k)}) \text{ so }
  (A^{(k)})^t A^{(k)}  (v-v^{(k)}) \sim (A^{(k)})^t  \big(\mM(v) - \mM(v^{(k)}) \big).
\end{displaymath}
We thus define the iterative algorithm
\begin{align}
\label{eq:Newton}
(A^{(k)})^t A^{(k)}  (v^{(k+1)}-v^{(k)}) = (A^{(k)})^t  \big(\mM(v) - \mM(v^{(k)}) \big).
\end{align}
Provided that the initial guess $v_0$ is sufficiently close to $v$ and that $A^tA$ is invertible, then it is a classical result that $(A^{(k)})^t A^{(k)}$ is also invertible and $v^{(k)}$ converges to $v$ as $k\to\infty$.

As we have seen in the preceding section, the  matrix $A^tA$ may not be invertible or may be ill conditioned;
for instance when $h_0^{''}(x_0)=0$ and $\alpha'_0(x_0)=0$, which imposes that $\psi_3=0$ in \eqref{eq:linearbis}.
To remedy this issue, we add a small penalty term and solve instead
\begin{align}
\label{eq:Newtonregul}
\Big( (A^{(k)})^t A^{(k)} + \lambda B^tB ) (v^{(k+1)}-v^{(k)}) = (A^{(k)})^t  \big(\mM(v) - \mM(v^{(k)}) \big),
\end{align}
with $B$ a linear operator that penalizes the curvature of $h$, for instance, $B$ can be a discretization of $h"(x)$; see \cite{engl} for references on regularization. We then choose $\lambda$ small in order not to affect the reconstruction of $h$ when the curvature does not vanish.

As a final remark, we note that nonlinear problems may be injective even when their linearizations are not (think of the map $x\mapsto x^3$ whose linearization is not invertible at $0$). However, when $h(x)=h_0$ and $\alpha(x)=\alpha_0$ are constant, then the available measurements ${\mathcal M}$ are independent of $h_0$.


\section{Numerical discretization of the graph formulation}
\label{sec:num}

\subsection{Description of the algorithm}

The upper boundary of the cloud can indeed be represented as a graph $(x,h(x))$ for a large class of clouds. However, the lateral and bottom boundaries would require different parameterizations. In order to avoid technical difficulties, we have considered clouds with the following (simplified and somewhat unrealistic) structure. We assume that the bottom of the cloud is flat  and given by the line segment between $(x_L,h_B)$ and $(x_R,h_B)$. The left part of the cloud is assumed vertical and given by the line segment between
$(x_L,h_B)$ and $(x_L,h_1)$. The right part of the cloud is also assumed vertical and given by the line segment between
$(x_R,h_B)$ and $(x_R,h_N)$. The upper part of the cloud is represented by $(x,h(x))$ for $x_L<x<x_R$ and discretized as a continuous piecewise linear function with nodes $(x_j,h_j)$ for $1\leq j\leq N$ such that $x_1=x_L$, $x_N=x_R$ and $x_j=x_L+\frac{j-1}{N-1}(x_R-x_L)$.

The outgoing radiant flux
(angularly-integrated radiance)
at each of the $N-1$ linear pieces of the upper surface is denoted by $\alpha_j$ for $1\leq j\leq N-1$. The fluxes 
on the left and right sides are denoted by $\alpha_L$ and $\alpha_R$, respectively. The function $\beta$ is discretized by a continuous piecewise linear function equal to $\beta_p$ at the angles $\pi\frac pP$ for some given value of $P$ equal to or comparable to $J$.

With this geometry for the cloud, it remains to construct a discretization of the forward map $\mM$. We assume measurements available for $\{\phi_j\}_{1\leq j\leq J}$ and discretize the upper surface by $X_n=nh$ for a discretization step $h$ and $n\in\Zm$. The value $u_{nj}$ is then the integral of the real solution $u(X,Z,\phi_j)$ generated by the discretized cloud averaged over the spatial interval $(nh,(n+1)h)$. We verify that $u_{nj}$ is non-zero for a finite number of values $n$.

The unknown parameters and (discretized) functions are therefore
\begin{displaymath}
  x_L,\,x_R,\,h_B,\, h_n,\,  \alpha_n,\, \alpha_L,\,\alpha_R,\,\beta_p,
\end{displaymath}
while the information is given by $u_{nj}$.

The reconstruction of $x_L$ and $x_R$ is performed as follows. We assume that $\phi_j=\frac\pi2$ for some $j$ and look at the support of $u_{nj}$. The smallest value $n_L$ such that $u_{n_Lj}\not=0$ and the largest value $n_R$ such that $u_{n_Rj}\not=0$ define $x_L=hn_L$ and $x_R=hn_R$.

The above procedure defines a continuous functional from $(h_B,\, h_n,\,  \alpha_n,\, \alpha_L,\,\alpha_R,\,\beta_p)$ to $u_{nj}$. That functional has been linearized as explained in the derivation of \eqref{eq:linear} and the corresponding nonlinear problem inverted as explained in section \ref{sec:nonlin}.

\subsection{Blocking of light}



For geometries with rapidly varying $h(x)$, some photons emitted from a point $x_0$ in a direction $\theta$ are such that the half line $0<s\mapsto x_0+s\theta$ intersects the boundary of the cloud before reaching the satellite line $(X\in\Rm,Z)$. In the numerical simulations, each half line $0<s\mapsto x_0+s\theta_j$ that starts from point $x_0$ is tested. We calculate the signed distances $(d_1,d_2)$ from every pair of consecutive points $(x_1,x_2)$ on the cloud, to the half line $0<s\mapsto x_0+s\theta_j$. The intersection happens when $d_1 d_2<0$. When it intersects the cloud before reaching the satellite line, then the effect of the radiation emitted from $(x_0,\theta_j)$ on the measurement $u(X,Z)$ is disregarded. The linear analysis of the preceding section shows that the reconstruction of $(\alpha,\beta,h)$ is possible so long as a sufficiently large number of directions $\theta_j$ reach the satellite line. This is confirmed by the numerical simulations presented below.

\subsection{Examples of numerical inversions}

We now present some numerical simulations that display the performance of the algorithm. 

\begin{figure}
\begin{center}
 \includegraphics[width=8cm]{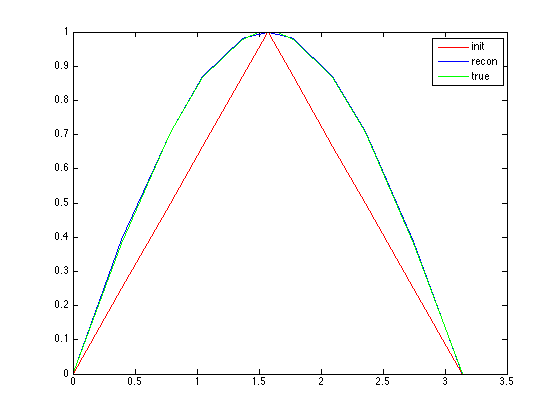}
 \end{center}
 \caption{True function (in green), reconstructed function (in blue) and initial guess (in red) for the angular radiation function $\beta(\phi)$.}
 \label{fig:beta}
\end{figure}

All reconstructions perform the reconstruction of $\alpha$, $\beta$, and the height $h$. In all algorithms, $\beta$ is given by the values presented in Fig. \ref{fig:beta} and represents a sine function. The initial guess and a typical reconstruction  for $\beta$ are also shown on that figure. 
For the rest of the section, we focus on the display of the reconstructions of $\alpha$ and $h$. In all simulations, the number of angular measurements is $J=9$. The spatial discretization is given by $N=51$ grid points.

In our examples, we used not only arbitrary functions for $\alpha$,
such as a constant or a step from zero to a non-zero value,
but also those that have more resemblance to real-world situations. Such cloud shape-dependent $\alpha$ functions are compiled as follows. Let $-\xi$
be direction above the horizon ($-\xi_z > 0$) from
which the sun sheds light on the cloud. 
Let $\nu(x)$ be the unit normal vector on the surface of the cloud that was defined earlier, and set $0 \leq \rho \leq 1$. Then,
\begin{equation}
\label{eq:solar_illum}
  \alpha(x)=\begin{cases}
    -\xi \cdot \nu(x), & \text{if $-\xi \cdot \nu(x) \geq \rho$,
                         and sun is not blocked};\\
    \rho, & \text{otherwise}.
  \end{cases}
\end{equation}\\
Thus  $\alpha(x) = \rho$ represents a self-shaded portion of the cloud's surface, where light emerges as a diffuse field,
while $\alpha(x) > \rho$ represents portions of the cloud's surface directly exposed to the sun.

In the figures below, we used three colors, blue, red and black to represent the true value, the reconstructed value and the initial guess respectively. In all cases the figures are 
displayed with $h(x)$ on the left and $\alpha(x)$ on the right.



{\bf Example 1}: $\xi=\frac{\pi}{6}$, meaning sun at 60$^\circ$ from zenith, $\rho=0.2$ (cloud with a lot of curvature)


\includegraphics[width=80mm, height=60mm]{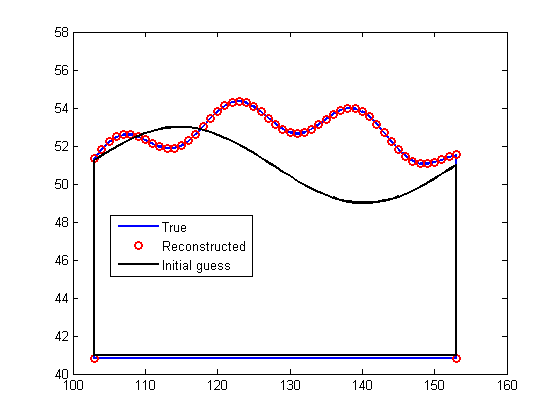}
\includegraphics[width=80mm, height=60mm]{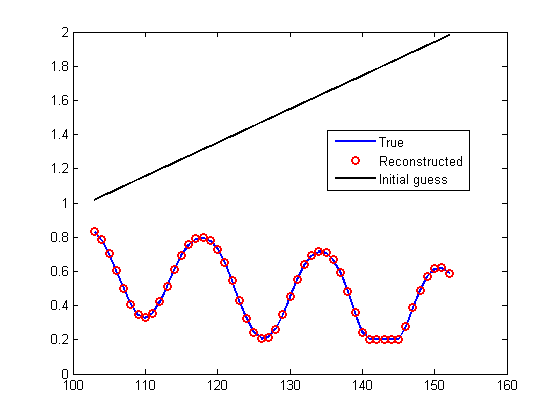}

{\bf Example 2}: $\xi=\frac{\pi}{2}$, meaning overhead sun, $\rho=0.2$ (cloud with a lot of curvature)

\includegraphics[width=80mm, height=60mm]{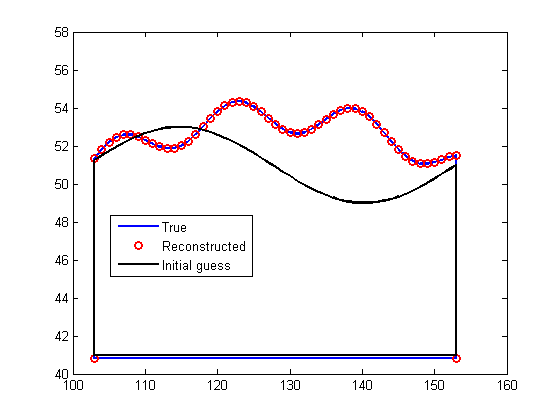}
\includegraphics[width=80mm, height=60mm]{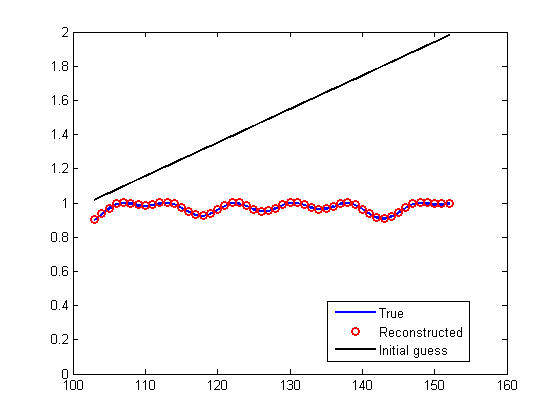}

{\bf Example 3}: $\xi=\frac{\pi}{2}$, $\rho=0.2$ (cloud with less curvature)

\includegraphics[width=80mm, height=60mm]{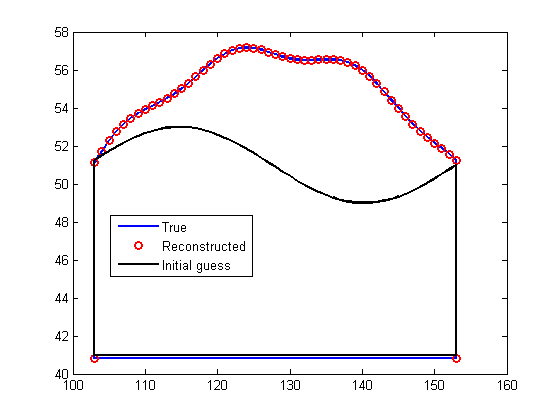}
\includegraphics[width=80mm, height=60mm]{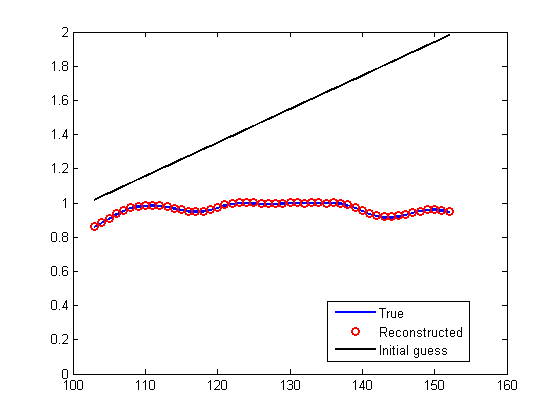}

{\bf Example 4}: In this case, $\alpha$ is chosen to be a constant (cloud with less curvature)

\includegraphics[width=80mm, height=60mm]{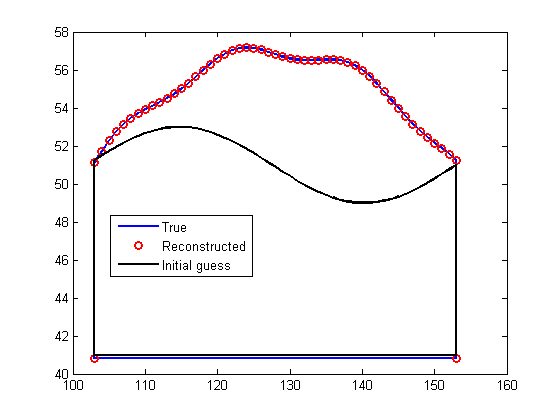}
\includegraphics[width=80mm, height=60mm]{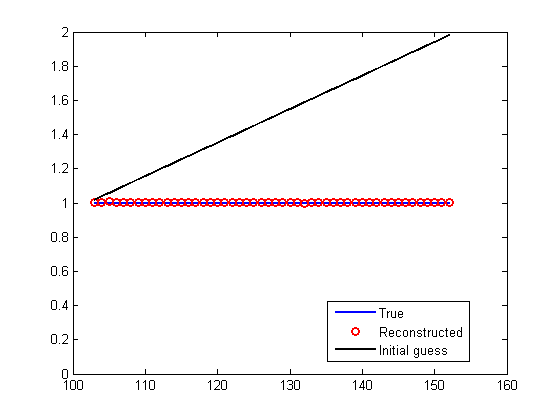}

{\bf Example 5}: $\alpha$ is chosen to be a step function (cloud with less curvature)

\includegraphics[width=80mm, height=60mm]{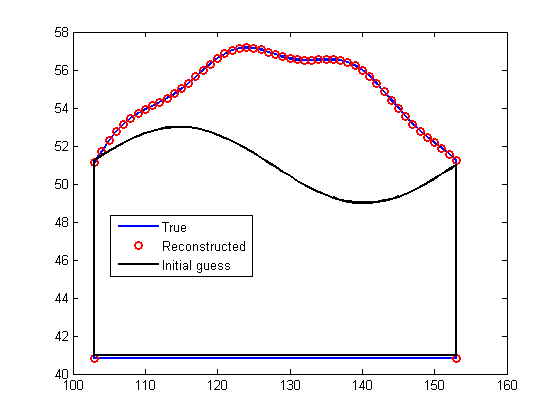}
\includegraphics[width=80mm, height=60mm]{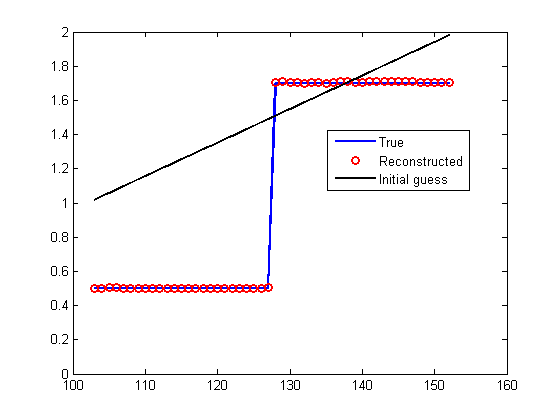}


In the five examples we presented above the algorithm works well in situations where $h(x)$ and $\alpha(x)$ are reasonably smooth, even for $\alpha(x)$ that has a jump discontinuity.

\section{Reconstruction of cloud speed}
\label{sec:speed}

Let us now assume that the cloud moves in the $x$ direction at constant speed $v_c$. We also assume that the satellite moves with speed $v_s$ along the $X$ axis and that $v_s>|v_c|$. (Both velocities are assumed to be far less than the speed of light.) When $v_c=0$, measurements are performed at $X=v_st$ for all values of $X$ or, equivalently, of $t$. This is the setting of the preceding section. When $v_c\not=0$, then the cloud is immobile in the reference $y=x-v_ct$. Define $Y=X-v_ct=(v_s-v_c)t = \frac{v_s-v_c}{v_s}X$. We define $\lambda=\frac{v_s-v_c}{v_s}$ an unknown parameter when $v_c$ is unknown (we assume $v_s$ known).

The inverse problem now consists of reconstructing $(\alpha,\beta,h,\lambda)$ from knowledge of $u$. Note that $\lambda$ is an additional scalar coefficient.

In the $(y,Y)$ variables, the problem is as before so that
\begin{align}\label{eq:yY}
  \lambda X-s(\lambda X,Z,\theta)\cos\phi &= x(\lambda X,Z,\theta), & Z-s(\lambda X,Z,\theta)\sin\phi &= h(x(\lambda X,Z,\theta)).
\end{align}
Let us denote $\tilde x(X,\phi)=x(\lambda X,\phi)$, dropping the dependency in $Z$ and again identifying $\theta$ and $\phi$. Introduce also $\tilde s(X,\phi)=s(\lambda X,\phi)$. Then all formulas in the preceding section with $x$ replaced by $\tilde x$ hold. In particular, with $x_0=x_0(\lambda X,\phi)$,
\begin{displaymath}
\begin{split}
 \delta u &= \beta_0\big(\phi-\psi_0(x_0)\big) \,\delta \alpha(x_0) + \alpha_0(x_0) \,\delta\beta(\phi-\psi_0(x_0)) \\ &
  + \alpha'_0(x_0)\beta_0(\phi-\psi_0(x_0))  \delta \tilde x(x_0) + \alpha(x_0) \beta'_0(\phi-\psi_0(x_0)) \dfrac{h_0"(x_0)}{1+(h_0'(x_0))^2} \delta \tilde x(x_0) \\
  &+     \beta'_0(\phi-\psi_0(x_0)) \dfrac{-\alpha_0(x_0) }{1+(h_0'(x_0))^2}  (\delta h)'(x_0)
\end{split}
\end{displaymath}

From \eqref{eq:yY}, we deduce that
\begin{displaymath}
  \delta h + h'_0 \delta \tilde x + \delta \tilde s \sin\phi =0 ,\quad \delta \tilde x + \delta \tilde s \cos\phi = X \delta\lambda.
\end{displaymath}
As a consequence, we have
\begin{displaymath}
    \delta \tilde x = \dfrac{\delta h + \tan \phi X \delta\lambda}{ \tan \phi - h'_0}
\end{displaymath}
We thus find, again with $x_0=x_0(\lambda X,\phi)$, that
\begin{align}
\label{eq:linearlambda}
  \delta u(X,Z,\phi) & = \beta_0\big(\phi-\psi_0(x_0)\big) \,\delta \alpha(x_0) + \alpha_0(x_0) \,\delta\beta(\phi-\psi_0(x_0)) \\ & + \dfrac{\alpha'_0(x_0)\beta_0(\phi-\psi_0(x_0)) + \alpha_0(x_0)  \psi'_0(x_0)\beta'_0(\phi-\psi_0(x_0))}{\tan \phi - h'_0(x_0)} \,\Big( \delta h(x_0) + X\tan \phi \delta\lambda \Big)\notag \\
  & +     \beta'_0(\phi-\psi_0(x_0)) \dfrac{-\alpha_0(x_0)}{1+(h_0'(x_0))^2}  (\delta h)'(x_0)\notag.
\end{align}
The reconstruction of $\delta\lambda$ is not always possible. Consider the case with $h'_0(x_0)=0$, i.e., the setting of a flat guess for the cloud. Then we observe, changing to $(X,\phi)$ to $(x_0,\phi)$ variables, that
\begin{displaymath}
  \delta u(x_0,\phi) =  \beta_0(\phi) \big(\delta\alpha(x_0) + \alpha_0'(x_0)X \delta\lambda\big) + \alpha_0 \delta\beta + \alpha_0' \dfrac{\beta_0(\phi)}{\tan\phi} \delta h - \alpha_0 \beta_0'(\phi) (\delta h)'
\end{displaymath}
Then as in the preceding section, we observe that $\delta\alpha(x_0) + \alpha_0'(x_0)X \delta\lambda$ can be reconstructed under appropriate assumptions on $\delta h$. Without prior information on $\delta\alpha$, it is impossible to separate $\delta\alpha$ from $\delta \lambda$.

In the simplified case where $h'_0\not=0$ and $h_0^{"}=0$ so that $\phi_0'=0$, and considering only the part involving $\delta \alpha$ and $\delta\lambda$, we find
\begin{displaymath}
    \beta_0(\phi-\psi_0(x_0)) \Big( \delta\alpha + \alpha'_0(x_0) \dfrac{X\tan\phi}{\tan\phi-h'_0(x_0)} \delta\lambda\Big).
\end{displaymath}
The functions $(1,\frac{\tan\phi}{\tan\phi-h'_0(x_0)})$ are then linearly independent as soon as $h'_0(x_0)\not=0$. In such a setting, $\delta\lambda$ and $\delta\alpha$ can be separately reconstructed, as well as $\delta h$ and $\delta \beta$ under conditions on $\delta h$ similar to those obtained in earlier sections.

This shows that in spite of what looks like redundant measurements ($u(X,\phi)$ known for all $X$ and $\phi$), the reconstruction of the additional scalar coefficient $\delta\lambda$ is not guaranteed in general. The nonlinear problem for $(\alpha,\beta,h,\lambda)$ is then handled as described earlier. That all coefficients in $(\alpha,\beta,h,\lambda)$ cannot be reconstructed has been confirmed in numerical experiments, where we were not able to obtain any converging algorithm. The simulations, carried out for values of $\lambda$ in the range 0.7 to 0.9, are not presented here.


Simultaneous reconstruction of the cloud's shape and speed requires a different modeling of the optical parameters of the cloud. This is the subject of ongoing research.

To motivate this future work, we note that horizontal cloud velocity is successfully retrieved at the same time as cloud (top) height using multi-angle imagery with 275~m resolution from the previously mentioned MISR instrument on NASA's Terra satellite \cite{HorvathDavies2001JTech,HorvathDavies2001GRL}.  In that case, $\lambda$ is very nearly unity since $v_s$ is an orbital velocity ($\sim$7~km/s) while $v_c$ is $\sim$7~m/s.  Apart from high-accuracy geolocation, the main requirement is to locate the same cloud or cloud feature as seen by different MISR cameras using standard feature matching methods \cite{MullerEtAl2002TGRS}, as is done routinely for cloud height estimation by stereography \cite{MoroneyEtAl2002TGRS}.  A minimum of three cameras is required to unravel cloud height and the along-track component of the wind.

\section{Reconstruction of more general geometries}
\label{sec:geom}

\subsection{Cloud boundaries in polar coordinates}

The constructions presented above assume that the cloud surface is represented by a graph $(x,h(x))$. Clouds with a closed boundary thus cannot be reconstructed by such a method unless several charts are considered. Several modifications of the theory allow us to reconstruct clouds with closed boundary. The simplest method is arguably to assume that the boundary can be represented as a graph in polar coordinates $[0,2\pi)\ni\theta\mapsto r(\theta)\in\Rm^2$. The measurements are now performed on a circle of radius $R$ that shares the same origin with the cloud, where $R> max(r)$ is fixed and known. Again, we assume that we receive measurements for $\{\phi_j\}_{1\leq j\leq J}$ and that $J\geq3$.  For an Earth science application, one can envision an aircraft carrying an imager that circumnavigates an opaque cloud of interest; this ``cloud'' could equally well be an opaque aerosol plume emanating from a powerful source (e.g., ash from an erupting volcano, smoke from a massive wildfire).

Let $u(\theta,r(\theta),\phi)$ be the density of photons as a function of position $(\theta,r(\theta))$ and angle $\phi$. Assuming the same conditions as in previous section, we have
\begin{displaymath}
u(\theta,r(\theta),\phi)=u(\Theta(\theta,r(\theta),\phi),R,\phi),
\end{displaymath}
where $(\Theta(\theta,r(\theta),\phi),R)$ is some position on the circle where we receive measurements.

Let $t(\theta,r(\theta))$ be the slope of the tangent line at a point $(\theta, r(\theta))$ of the surface, and let
\begin{equation} \label{eq:normal}
\nu(\theta,r(\theta)) = \dfrac1{\sqrt{1+\Big(t\big(\theta,r(\theta)\big)\Big)^2}} \,\Big(-t\big(\theta,r(\theta)\big),1\Big)
\end{equation}
be the outward unit normal vector to the surface. Then as before, we assume that
\begin{displaymath}
u(\theta,r(\theta),\phi)=\alpha(\theta)\beta\Big(\phi-\arctan\big(t(\theta,r(\theta),\phi)\big)\Big),
\end{displaymath}
for two unknown functions $\alpha(\theta)$ and $\beta(\phi)$. From above we thus observe that,
\begin{displaymath}
u(\Theta,R,\phi)=\alpha\big(\theta(\Theta,R,\phi)\big)\beta\Big(\phi-\arctan\big(t(\theta(\Theta,R,\phi))\big)\Big),
\end{displaymath}
where we still assume that $\beta(\phi)$ is normalized.

In this setting, we thus have measurement operator
\begin{align}
\label{eq:measop}
{\mathcal K}: \, (\alpha,\beta,r) \mapsto {\mathcal K}(\alpha,\beta,r)  = u(\Theta,R,\phi).
\end{align}
The inverse problem consists of reconstructing $(\alpha,\beta,r)$ from knowledge of ${\mathcal K}(\alpha,\beta,r)  = u(\Theta,R,\phi)$.

\subsection{Description of the algorithm}

Let us represent the cloud boundary by a graph $(\theta,r(\theta))$ for $\theta \in[0,2\pi)$, which we assume to be piecewise linear. For $1 \leq j \leq N+1$, we have $r(\theta_1)=r(\theta_{N+1})$ and $\theta_j=\theta_1+\frac{j-1}{N}2\pi$.

The radiance at each of the $N$ linear pieces of the cloud is denoted by $\alpha_j$ for $1 \leq j \leq N$. The function $\beta$ is discretized in the same way as previously.

With this geometry for the cloud, it remains to construct a discretization of the forward map $\mathcal K$. We assume measurements available for $\{\phi_j\}_{1 \leq j \leq J}$ and discretize the measurement circle of radius $R$ by $\Theta_n=n\Delta\Theta$ for a discretization step $\Delta\Theta$ and number $n(R) \in \Zm$. We denote $n$ as $n(R)$ here because the bigger the radius the finer we have to discretize in order to get a reasonably smooth measurement. The value $u_{nj}$ is then the integral of the real solution $u(\Theta,R,\phi_j)$ generated by the discretized cloud averaged over the arc $(\Theta_n R,\Theta_{n+1} R)$. We verify that $u_{nj}$ is non-zero for a finite number of values $n$.

The above procedure defines a continuous functional from $(\theta_m, r_m, \alpha_m, \beta_p)$ to $u_{nj}$. The functional can be linearized similarly as explained in the derivation of \eqref{eq:linear} and the solution of the corresponding nonlinear inverse problem is similar to what has been explained in section \ref{sec:nonlin}.

\subsection{Numerical simulations}

We now show some numerical simulations generated from the model we just described. In our settings, we chose $(0,0)$ to be the origin of our cloud and of the disk of radius $R$ from which we receive measurements. In our examples, we used not only arbitrary functions for $\alpha$, but also the more realsitic one in \eqref{eq:solar_illum} that mimics solar illumination.
In the following figures, we continue to use three colors, blue, red and black to represent the truth, the reconstruction and the initial guess, respectively. Our initial guess for $\alpha$ is always a constant function $\alpha(\theta)=1$. In all cases, the figures would be displayed with $r(\theta)$ on the left and $\alpha(\theta)$ on the right.

The function $\beta$ is still as in Fig.\ref{fig:beta}. The value of $N=201$ (so that $d\theta=2\pi/200$) and the number of measurement directions is $J=11$. The angles are given explicitly by the arccos value of $(\pm1,\pm\cos \frac\pi4,\pm\cos \frac\pi3,\pm\cos \frac\pi{2.3},\pm\cos \frac\pi{2.1},0)$.

{\bf Example 1}: A smooth reconstruction of $r(\theta)$ and $\alpha(\theta)$:

\includegraphics[width=80mm, height=60mm]{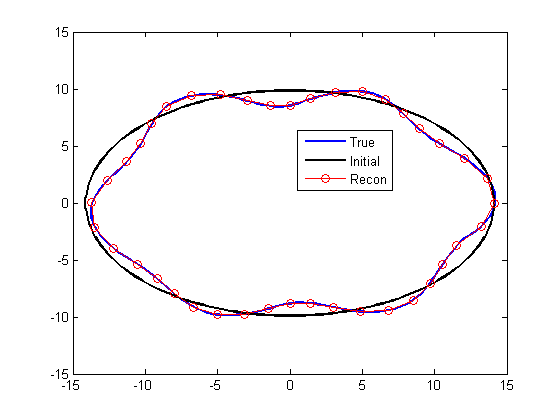}
\includegraphics[width=80mm, height=60mm]{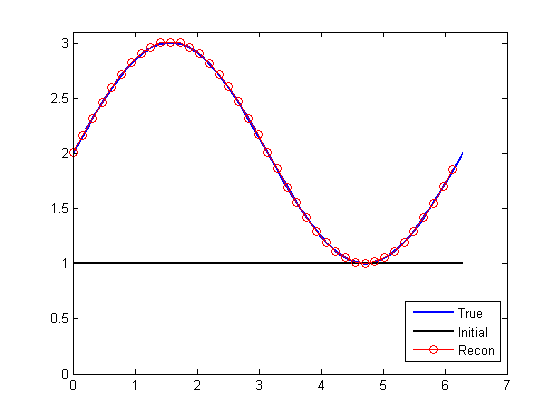}

In this example, we obtain excellent reconstructions of both $r(\theta)$ and $\alpha(\theta)$.

{\bf Example 2}: $\xi=\frac{\pi}{6}$, $\rho=0.2$ and $\alpha$ given in \eqref{eq:solar_illum}.

\includegraphics[width=80mm, height=60mm]{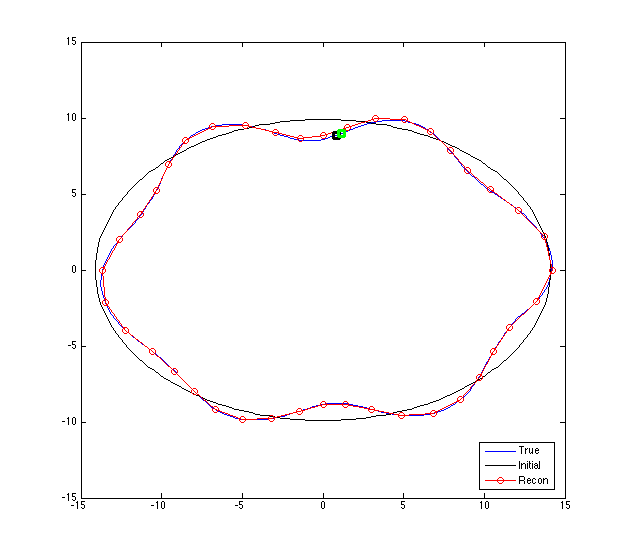}
\includegraphics[width=80mm, height=60mm]{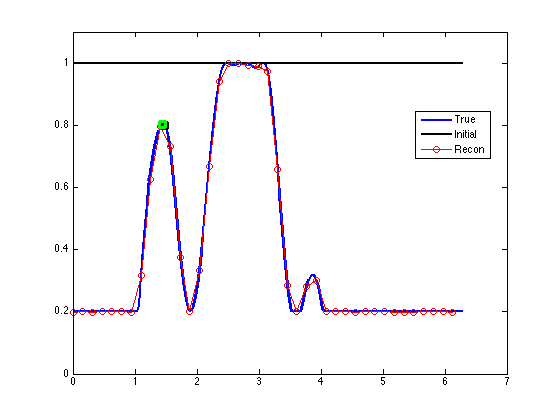}

In this example, the reconstructions for both $r(\theta)$ and $\alpha(\theta)$ are good except in the vicinity of the two points highlighted in green and black boxes in both figures (upper middle to the right). Note that we have mentioned earlier that the reconstruction would break down when $\alpha'=0$ and $h''=0$, where $(x,h(x))$ represents a Cartesian representation of the cloud surface in an appropriate system of coordinates. We could see clearly that some point between the highlighted points is precisely where $\alpha'=0$ and $h''=0$. The reconstruction would have been worse had we not added the small penalty term that was described earlier in section \ref{sec:nonlin}.  


{\bf  Example 3}: $\xi=\frac{5\pi}{6}$, $\rho=0.2$  and $\alpha$ given in \eqref{eq:solar_illum} (a well-reconstructed profile)

\includegraphics[width=80mm, height=60mm]{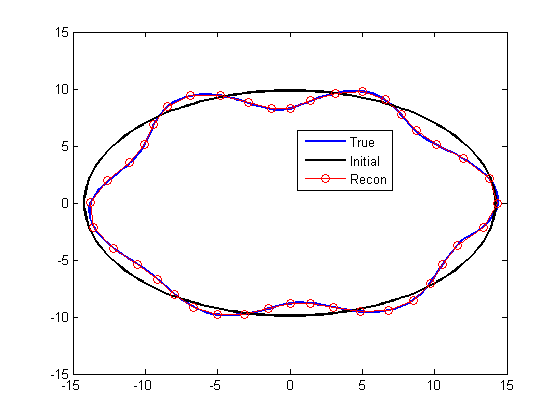}
\includegraphics[width=80mm, height=60mm]{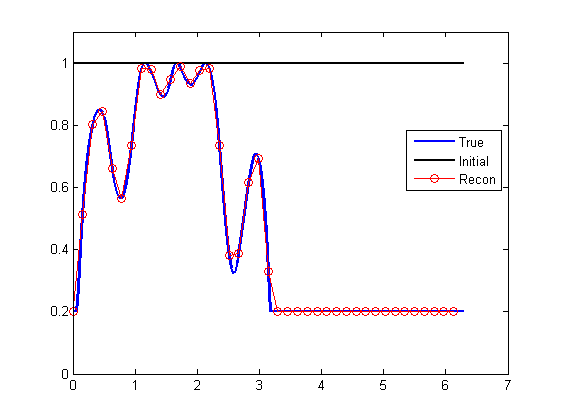}

In this configuration, all hypotheses necessary for a good reconstruction are met.


{\bf Example 4}: In this case, $\alpha$ is chosen to be a step function

\includegraphics[width=80mm, height=60mm]{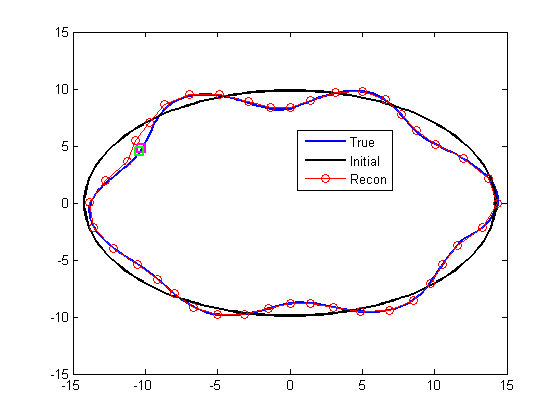}
\includegraphics[width=80mm, height=60mm]{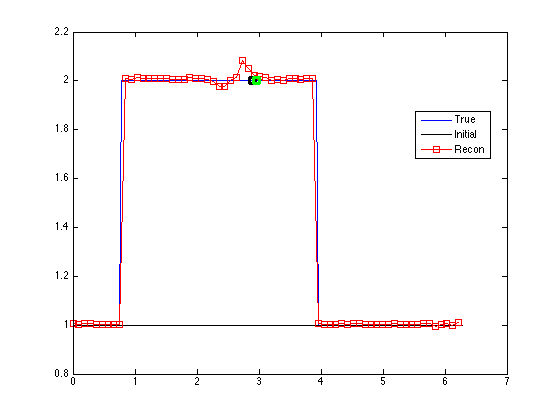}

In this example, we gave more curvature to the cloud and again the reconstructions were good for most of the parts except in the vicinities of the points highlighted in green and black boxes, which were highlighted in both figures (lower to the right). We can see clearly that the reconstruction broke down for the same reason as in Example 2.


%


\section{Conclusions and outlook}
\label{sec:conclu}

This paper presents reconstruction procedures for the geometry of a cloud from airborne or satellite observations at $J \ge 3$ viewing angles and at quite high spatial resolution (on the order of 100s of pixels across the cloud). Rather than reconstructing the optical properties of cloud, which may vary in space, we focus here on a robust reconstruction of the geometry of the cloud without modeling light scattering.

In the simplified setting of a two-dimensional cloud with known speed we showed that the three one-dimensional functions $(\alpha,h,H)$ (or equivalently $(\alpha,h,\beta)$; see \eqref{eq:chbeta}) could all be reconstructed from satellite measurements when $J$ is sufficiently large in favorable situations where the linearized map is invertible. This was confirmed by several numerical solutions showing that our iterative procedure converged in many settings. However, there are situations where such a map is not invertible. Specifically, convergence problems arise when the cloud's boundary lacks curvature and/or the angular distribution of light emerging from the cloud is nearly isotropic.

More surprising, we found that the speed of the cloud, when unknown, could not always be uniquely recovered from the available measurements.  This finding immediately motivates further research since cloud-based wind speeds are routinely retrieved from multi-angle imaging sensors in space with far lower spatial resolution than assumed here (on the order of 10s of pixels across the cloud, or less).  On the other hand, no attempt is made yet with these sensors to retrieve cloud shape, only their height in some spatially-averaged sense.

We expect many results to hold in a three-dimensional environment. The main new feature is that the now three-dimensional normal vector $\nu(x)$ is no longer in the plane given by the satellite trajectory and measured directions $\theta_j$, which would slightly modify the form of the linearized operator mapping $\delta v$ to $\delta u$.

The main result of the paper is to show that geometric features of the cloud could be reconstructed from satellite measurements without having to estimate the radiative transfer parameters inside a cloud. In some configurations, for instance when the cloud speed is unknown, this procedure may have too strong limitations. In such a case, a parameterization of the optical properties of the cloud becomes necessary. The resulting inverse problem for $(\alpha,h,H)$ and such optical properties is the object of current research.

\section*{Acknowledgment}
This work was partially funded by NSF grant DMS-1408867 and two NASA grants, one from SMD/ESD/Radiation Sciences (programs managed by Hal Maring and Lucia Tsaoussi) and one from ESTO/AIST (managed by Michael Seablom).

%
%
%


\end{document}